\documentclass{amsart}

\usepackage{amssymb,amsmath,amsthm,amscd,float} 
\allowdisplaybreaks

\usepackage{booktabs, topcapt}

\usepackage{hyperref}

\hypersetup{
  colorlinks   = true, 
  urlcolor     = blue, 
  linkcolor    = blue, 
  citecolor   = red 
}

\usepackage[sorted, msc-links]{amsrefs}

\usepackage{mathtools}
\newtheorem{thm}{Theorem}
\newtheorem{prop}{Proposition}
\newtheorem{lem}{Lemma}
\newtheorem{rem}{Remark}

\newtheorem{exa}{Example}
\newtheorem{prob}{Problem}
\theoremstyle{defn}

\newcommand\M{\mathcal M}
\def\O{\mathcal O}

\def\cR{\mathcal R}

\def\l{\lambda}

\def\Z{\mathbb Z}
\def\Q{\mathbb Q}
\def\C{\mathbb C}
\def\R{\mathbb R}
\def\P{\mathbb P}

\def\H{\mathcal H}
\def\M{\mathcal M}

\def\l{\lambda}

\def\<{\langle}
\def\>{\rangle}
\def\X{\mathcal X}

\def\D{\Delta}

\def\X{\mathcal X}

\def\ch{\mbox{\rm char } }
\def\J{\mathcal J}

\def\A{\mathbb A}
\def\h{\mathfrak h}

\newcommand\p{\mathfrak p}

\def\WP{\mathbb{WP}_{w}^6}

\begin{document}

\title{On  hyperelliptic  curves of  genus 3}

\author{L. Beshaj}
\address{Department of Mathematical Sciences \\ United States Military Academy\\ West Point}
\email{lubjana.beshaj@usma.edu}

\author{M. Polak}
\address{Department of Computer Science \\ Rochester Institute of Technology}
\email{mkp@cs.rit.edu}


\date{}                                           


 

\begin{abstract} 
We study the moduli space of genus 3 hyperelliptic curves   via the weighted projective space of binary octavics. This enables us to create a database of all genus 3 hyperelliptic curves defined over $\Q$, of weighted moduli height $\h =1$. 
\end{abstract}

\maketitle

\tableofcontents

\section{Introduction}

Genus 3 curves have been interesting objects in classical mathematics for  many reasons. They are the lowest genus curves which are not hyperelliptic. Their generic equation is a ternary quartic given explicitly as in \cite{generic-3}. The lowest genus curve, defined over a field of characteristic zero,  which reaches the maximum size for the automorphism group is e genus three curve, namely the well celebrated Klein's quartic which has automorphism group of order 168.  

The moduli space of genus three curves $\M_3$ is a coarse moduli space of dimension 6. The hyperelliptic locus $\H_3$  in $\M_3$ has dimension 5.  It is precisely this fact that makes genus 3 curves very attractive in certain cryptographical applications; see \cite{frey-sh} for details. The existence of isogenies between genus 3 hyperelliptic Jacobians and genus 3 non-hyperelliptic Jacobians are the focus of much research lately; see \cite{c-m-sh} where   fibrations of such surfaces are studied via the Satake polynomial as in \cite{Satake}.

This paper started with the idea of creating a database of genus 3 hyperelliptic curves similar to the database for genus 2 curves in \cite{data}.  There are a couple of approaches to this.  In \cite{nato-1} and \cite{nato-2} we created such database for genus 2 curves with small coefficients by making use of the notion of height in the projective space.  Such techniques could be used for any genus $g>2$ hyperelliptic curves providing that one can handle invariants of binary forms. However, a new important point of view has become available after Shaska in \cite{sh-h} introduced the concept of \textit{weighted moduli height} for weighted moduli spaces. This makes it possible to created a database with small invariants instead of small coefficients. 

There are several benefits in using the weighted projective moduli space instead of the regular moduli space.  First, using the weighted moduli space we avoid special cases that we have to address during the compactification of $\H_3$. Another benefit is computational: one has to store powers of invariants when dealing with $\H_3$, but this is not necessary when dealing with the weighted moduli space.  The case of $g=2$ treated in \cite{m-sh} illustrates how things become much easier when using the weighted moduli space instead of the regular moduli space.  


The main idea of this paper is to use the approach from \cite{m-sh} for binary sextics and to study the weighted moduli space of binary octavics, including their twists and the weighted moduli height. We will follow the definition of the weighted moduli height as in \cite{m-sh}.  There are a few things that are slightly different for binary octavics from the binary sextics, as we will see in the coming sections. 

It is important to notice that the set of generators of the ring of invariants of binary octavics (or the so-called \textit{Shioda invariants}), namely $J_2$, $J_3$, $J_4$, $J_5$, $J_6$, $J_7$, $J_8$ satisfy five syzygies among them \cite{shioda} which were determined by Shaska in \cite{gen-3} as a single degree 5 equation.  Hence, a tuple $(J_2, \dots , J_8$ doesn't necessary correspond  genus 3 hyperelliptic curve.  There are two conditions that need to be checked here.  First, if such tuple satisfies the equation $F(J_2, \dots , J_8)=0$ given in \cite{gen-3} and second if the corresponding discriminant $J_{14}$ is non-zero.  

This paper is organized as follows. In section two we define the invariants of the binary octavics in terms of transvections and in terms of root differences, see \cite{shioda,  gen-3} for more details.   As we will see, two genus three curves $\X$ and $\X^\prime$ are isomorphic   over an algebraically closed field  $k$ if and only if there exists some $\l \in k\setminus \{ 0\}$ such that 
\[ J_i (\X) = \l^i J_i(\X^\prime), \textit{   for   }  \,\,  i=2, \dots , 7,   \]  
and the invariants  $J_2, \dots, J_8$ satisfy \cite[Eq.~(14)]{gen-3}. Hence, isomorphic curves correspond to elements of the projective space $\P_{(2,3,4,5,6,7,8)}$ such that the coordinates satisfy the two conditions mentioned above.  Therefore, enumerating points on the moduli space can be done by enumerating points on a variety defined inside the weighted projective space.  Hence, we can built a database of isomorphism classes of genus three hyperelliptic curves using $[J_2, \dots, J_8]$ as a point in the weighted projective space $\WP (\Q)$, with weight  $w = (2, \dots, 8)$,  satisfying  \cite[Eq.~(14)]{gen-3}.  

In section three, we describe the weighted moduli space of binary octavics following closely the notation in \cite{sh-h}.  Then we define the weighted moduli height which will be used in the later section to construct the database.  Using this height for $\WP (\Q)$,  $w = (2, \dots, 8)$, we are able to determine a unique tuple $(J_2, \dots, J_8)$ for a point $ \p \in \WP (\Q)$ and compile a database of genus three hyperelliptic curves.

In the last section we describe in details how we can build a database of genus three hyperelliptic curves.

\medskip

\noindent \textbf{Aknowledgments:} We want to thank Prof. Shaska for sharing his Maple package on genus 3 hyperelliptic curves, including the implementation in Maple of the equation \cite[Eq.~(14)]{gen-3}.  We also want to thank Scott Guest for sharing his Sagemath package on weighted projective spaces.

\newpage 
\section{Preliminaries}
Let $k$ be  an algebraically closed field.   A binary  form  of degree $d$ is  a homogeneous  polynomial $f(X,Y)$ of  degree $d$ in two  variables over $k$.  Let $V_d$ be the $k$-vector space of binary  forms of degree $d$.  The group $GL_2(k)$ of  invertible $2  \times 2$  matrices over  $k$  acts on  $V_d$  by coordinate  change. Many  problems  in algebra  involve properties  of binary forms  which are invariant under these  coordinate changes. In particular, any  hyperelliptic genus $g$ curve over  $k$ has a  projective equation of the form $Z^2Y^{2g} = f(X,Y)$,  where $f$ is  a binary form  of degree  $d=2g+2$  and  non-zero discriminant. Two  such curves  are isomorphic  if and  only if  the corresponding  binary forms are conjugate under $GL_2(k)$. Therefore the moduli space $\mathcal H_g$ of hyperelliptic genus $g$ curves  is the affine  variety whose  coordinate ring  is the  ring of $GL_2(k)$-invariants in the coordinate ring  of the set of elements of $V_d$ with non-zero discriminant.  
Throughout this section $\ch (k) \neq 2, 3, 5, 7$.

\subsection{Covariants and invariants of binary octavics}
We will use the symbolic method of classical theory to construct covariants of binary octavics.  Our standard reference for this section is  \cite{shioda} or \cite{gen-3}. 
Let $f$   denote a binary octavic then  $f (X, Y )$ is given by the formula:

\[f (X, Y ) =\sum_{i=0}^8a_iX^i Y^{8-i}.\]

For any two binary forms $f, g$ of degree of degree $n$ and $m$ respectively the level $r$ transvection is 
 \[(f, g)^r=\dfrac{(m-r)!(n-r)!}{n!m!}\sum_{k=0}^r(-1)^k{{r}\choose{k}}\cdot \dfrac{\partial^rf}{\partial X^{r-k}\partial Y^k}\cdot\dfrac{\partial^rg}{\partial X^{r}\partial Y^{r-k}}\]

As in Eq.~(8) in \cite{gen-3}  we define the following
covariants:
\begin{equation}
\begin{split}
&g=(f,f)^4, \quad k=(f, f )^6, \quad h=(k,k)^2, \quad m=(f,k)^4, \\
&   n=(f,h)^4, \quad p=(g,k)^4, \quad q=(g, h)^4.\\
\end{split}
\end{equation}

\noindent Then, the following 
\begin{equation}\label{def_J}
\begin{aligned}
&     \J_2=  (f,f)^8,                    &   \qquad   &  \J_3=  (f,g )^8, \qquad  \J_4=  (k,k)^4,                      &    \qquad      &   \J_5=  (m,k)^4,  \\
&   \J_6 =   (k,h )^4,              &    \qquad      &      \J_7=  (m,h )^4,  \qquad \J_8=   (p,h)^4,        &    \qquad       &   \J_9=    (n,h)^4, \\
&  \J_{10}=   (q,h)^4      &  \qquad    &       \\
\end{aligned}
\end{equation}
are $SL_2(k)$- invariants.     There is another invariant $\J_{14}$ which is the discriminant of the binary form $f(X, Y)$.  For the purposes of this paper we will denote $\D = \J_{14}$.

\begin{thm}
The graded ring $\cR_8$  of invariants of binary octavics is generated by elements $\J_2, \dots , \J_{10}$. 
\end{thm}
 
See \cite[Thm.~3]{shioda} for a computational proof and \cite[Theorem~5]{gen-3} for a proof via the Reynold's operator.    
Moreover, the ring of invariants  $\cR_8$ of binary octavics is finitely generated as a module over $k[\J_2, \dots,\J_7]$; see \cite{shioda} or \cite{gen-3}.

In this paper we will follow closely the computations on \cite{gen-3}. For reasons that will be explained in the following section  we prefer that  the definition of $[\J_2: \J_4: \dots : \J_8]$ constitutes an  integer tuple.  Thus, we modify the definitions as follows. 
\begin{equation}\label{def_J}
\begin{aligned}
&     J_2= 2^2 \cdot 5 \cdot 7 \cdot (f,f)^8,                    &   \qquad   &  J_3= \frac 1 3 \cdot  2^4 \cdot 5^2 \cdot 7^3 \cdot (f,g )^8, \\
&  J_4= 2^9 \cdot 3 \cdot 7^4 \cdot (k,k)^4,                      &    \qquad      &   J_5= 2^9 \cdot 5 \cdot 7^5 \cdot (m,k)^4,  \\
&   J_6 = 2^{14} \cdot 3^2 \cdot 7^6 \cdot (k,h )^4,              &    \qquad      &      J_7= 2^{14} \cdot 3 \cdot 5 \cdot 7^7 \cdot (m,h )^4,  \\
& J_8= 2^{17} \cdot 3 \cdot 5^2 \cdot 7^9 \cdot   (p,h)^4,        &    \qquad       &   J_9= 2^{19} \cdot 3^2 \cdot 5 \cdot 7^9 \cdot   (n,h)^4, \\
&  J_{10}=   2^{22} \cdot 3^2 \cdot 5^2 \cdot 7^{11} (q,h)^4      &  \qquad    &       \\
\end{aligned}
\end{equation}

In \cite[Theorem 6]{gen-3} it was shown that    invariants  $J_2, \dots , J_8$ satisfy  the following equation
%
\begin{equation}\label{shaska}
 J_8^5 + \frac { I_8} {3^4 \cdot 5^3 } J_8^4 + 2 \cdot \frac { I_{16} } {3^8\cdot 5^6}   J_8^3 + \frac {I_{24}} {2 \cdot 3^{12} \cdot 5^6}  J_8^2  + \frac { I_{32}} {3^{16} \cdot 5^{10} }  J_8 + \frac { I_{40}} {2^2 \cdot 3^{20} \cdot 5^{12}}  =0, 
 \end{equation}
where  $I_8, I_{16}, I_{24}, I_{32}$ and $I_{40}$ are displayed in \cite[Appendix]{gen-3}.

Thus, in our efforts to create a list of genus 3 hyperelliptic curves we will create a list of ordered tuples $\left( J_2, \dots , J_8 \right)$ such that its coordinates satisfy Eq.~\eqref{shaska}.  When do two such tuples give the same curve?  We have the following:

\begin{prop} Two genus 3 hyperelliptic curves $\X$ and $\X^\prime$
%
%
are isomorphic over $k$ if and only if there exists some $\l \in k\setminus \{ 0\}$ such that 
\[ J_i (\X) = \l^i J_i(\X^\prime), \textit{   for   }  \,\,  i=2, \dots , 7,   \]  
and $J_2, \dots, J_8$ satisfy the \cite[Eq.~(14)]{gen-3}.
Moreover, the isomorphism is given by  $f=f^M$   where $M  \in GL_2(k) $ and $\l = \left( \det M \right)^4$.
\end{prop}

See   \cite[Theorem~7]{gen-3} for details of the proof. 
In the next section we will explain how the tuples $\left( J_2, \dots , J_8 \right)$ are the points in the weighted  projective space.  


There is no known rational set of generators for the graded ring $\R_8$.    In \cite{gen-3} Shaska has defined a set of absolute invariants $t_1, \dots , t_6$ which seem to work very well for curves of relatively small naive height, but as expected they are not defined everywhere.  

\subsection{Invariants from root differences} 

In 1986  Tsuyumine, while studying the Siegel modular forms of degree 3,  gave a proof that the graded ring $S(2, 8)$ of invariants of binary
octavics is generated by invariants expressed in terms of root differences.
 Let $\xi_1,...,\xi_8$ denote the roots of $\sum {8\choose i} a_ix^{8-i}$ and $(ij)$ denotes the roots difference $\xi_i-\xi_j$. We define the following as in \cite[pg. 772]{ts-1}.
\begin{equation}
\begin{split}
I_2 & = \displaystyle\sum (13)(14)(23)(24)(57)(58)(67)(68), \\
I_3 & = \displaystyle\sum (12)^2(34)^2(56)^2(78)^2(13)(24)(57)(68),\\
I_4 & = \displaystyle\sum (12)^4(345,678)^2,\\
I_5 & = \displaystyle\sum (12)^4(345,678)^2(15)(26)(37)(48),\\
I_6 & = \displaystyle\sum (1234,5678)^2,\\
I_7 & = \displaystyle\sum (1234,5678)^2(15)(26)(37)(48),\\  
I_8 & = \displaystyle\sum (1234,5678)^2(15)(16)(25)(26)(37)(38)(47)(48),\\
I_9 & = \displaystyle\sum (1234,5678)^2(15)(16)(17)(26)(27)(28)(35)(36)(38)(45)(46)(48),\\
I_{10} & = \displaystyle\sum (1234,5678)^2(15)^2(26)^2(37)^2(48)^2(16)(17)(25)(28)(35)(38)(46)(47),
\end{split}
\end{equation}
where $(1234,5678)$ denotes $(12)(13)(14)(23)(24)(34)\times(56)(57)(58)(67)(68)(78)$ and $(345,678)$ denotes
$(34)(35)(45)(67)(68)(78)$.
Then we have the following.

\begin{prop}[Tsuyumine]
The graded ring $\cR_8$ of invariants of binary octavics is generated by $I_2, I_3, I_4, I_5, I_6, I_7, I_8, I_9, I_{10}$. 
\end{prop}

In \cite{ts-1} all degree 3 Siegel modular forms are expressed in terms of invariants $I_2, \dots , I_{10}$ and in terms of thetanulls.  
In \cite{p-sh-w} we display all thetanulls for genus 3 hyperelliptic curves.  Moreover, the thetanull constraints that define the loci for each automorphism group are calculated.

As far as we are aware, the explicit formulas  converting $I_2, \dots , I_{10}$ to $J_2, \dots , J_{10}$ and vice-versa have not appeared in print.

\begin{prob}
Express $I_2, \dots , I_{10}$ in terms of  $J_2, \dots , J_{10}$ and vice-versa
\end{prob}

One can express the above sets of invariants in terms of the Siegel modular forms, but that is outside the scope of this paper.  
In the next section, we will describe the weighted projective space and define the weighted moduli height. 
\section{Weighted moduli space of binary octavics}
 Let $k$ be a field of characteristic zero and  $q_0, \dots , q_n$ a fixed tuple of positive integers called \textbf{weights}. 
 Consider the action of $k^\star$ on $\A^{n+1}$ as follows
\[ \lambda \star (x_0, \dots , x_n) = \left( \l^{q_0} x_0, \dots , \l^{q_n} x_n   \right) \]
for $\l\in k^\ast$. 
The quotient of this action is called a \textbf{weighted projective space} and denoted by   $\P^n_{(q_0, \dots , q_n)}$. 
It is the projective variety $Proj \left( k [x_0,...,x_n] \right)$ associated to the graded ring $k [x_0, \dots ,x_n]$ where the variable $x_i$ has degree $q_i$ for $i=0, \dots , n$.   Next we focus our attention to the weighted projective space of binary octavics.
\subsection{The weighted moduli space of binary octavics}

From above we know that the invariant ring of binary octavics is generated by invariants $J_2, \dots , J_8$.
Since $J_i$, for $i=2, \dots , 7$ are all homogenous polynomials of degree $i$, we take the set of weights $w =(2,3,4,5,6,7, 8)$ and considered the weighted projective space $\WP (k)$. 

Invariants of octavics define a point in a weighted projective space $\p = [J_2 :J_3: J_4 :J_5: J_6 : J_7 : J_8 ] \in  \WP $. However, not every point in $\WP$ correspond to  a   genus 3 hyperelliptic curve.   

\begin{prop}\label{prop-1} 
Let $\p \in \WP$, where $w =(2,3,4,5,6,7, 8)$.  Then $\p$ correspond to the isomorphism class of a genus 3 hyperelliptic curve if and only if    $\p \in \WP (k) \setminus \{\D = 0\}$ and its coordinates satisfy Eq.~\eqref{shaska}.  
\end{prop}

Hence, in our goal of creating a list of all points in $\WP$ which correspond to genus 3 hyperelliptic curves we will have to through away all the points left out by the Prop.~\ref{prop-1}.  But first, let's determine a way of ordering such points in $\WP$. 

\subsection{The height of $\WP$}
Let $ K \subset k$ be a number field and $\O_K$ its ring of  integers. For a point $\p$ such that  $\p = [J_2 : J_3 : J_4: J_5: J_6 : J_7 : J_8 ] \in \WP (\O_K )$, we say that  $\p$ is defined over $\O_K$.  We call the tuple $\p$ a \textbf{minimal tuple} or a \textbf{normalized  weighted moduli point} if there is no prime $p \in \O_K$  such that $p^{ i }| J_{i}$, for  $i = 2,3,4,5,6,7, 8$. Two minimal tuples $\p$ and $\p^\prime$   are called \textbf{twists} of each other if there exists $\lambda \in k$ such that 
\[ \lambda \star \p = \p^\prime.\]

Let $\p = [J_2 : J_3 : J_4 : J_5: J_6 : J_7 : J_8 ]$ be any tuple in $\WP (K)$. We define the \textbf{weighted moduli height} of $\p$ (or simply \textbf{the height}) to be 

\begin{equation}\label{height} 
\h (\p)= \dfrac{\max\{|J_i|^{\frac{1}{i}}\}  }{\displaystyle \prod_{p^i|J_i}p},
\end{equation}
where the product is taken over all the primes $p \in\O_K$. In \cite{m-sh} was proved that this is a well-defined height in a
weighted projective space. The weighted moduli height of a minimal tuple is  simply $\h (\p)= \max\{|J_i|^{\frac{1}{i}}\}$. It is easy to verify that:

\begin{lem} Let   $\p=[J_2, \dots , J_7, J_8] \in \WP$. Then the following hold:

i)   $\p$ is  a minimal tuple if and only if it has minimal height. 

ii) If $\p$ is a minimal tuple and it has    a twist with minimal height, then there exists a square free integer $d$ such that  $d^i | J_{i}$ for $i = 2, \dots , 8$.
\end{lem}


A point $\p=[J_2, \dots , J_8] \in \WP$ is called an  \textbf{absolute minimal tuple} if it has  the smallest height among all the twists.

\begin{lem}
A minimal tuple  $\p=[J_2, \dots , J_8] \in \WP (\Q)$ is  an absolute minimum tuple   if and  only if there is no 
$\lambda \in \C\setminus\{0\}$ such that $\l^i \in \Z$ and  $\l^i | J_i$ , for all $i =2, \dots , 8$ such that $J_i\neq0$.

\end{lem}

\proof
Let $\p$ be an absolute minimal tuple. Suppose that there exist $\lambda \in \C\setminus\{0\}$ such that $\l^i \in \Z$ and  $\l^i | J_i$, for all $i =2, \dots , 8$ such that $J_i\neq0$.
Hence, 
\[\p^\prime=\left[ \frac{J_2}{\l^2},\dots,\frac{J_8}{\l^8}\right]=\left[J_2^{\prime}, \dots , J_8^{\prime}\right] \in \WP (\Z) \]
is a twist of $\p$ such that $J_i = \l^i J_i^\prime$, for all $i=2, \dots, 8$. Therefore, $J_i\geq J_i^\prime$ for all $i=2, \dots ,  8$ and 
\[\h(\p)=\max\{|J_i|^{\frac{1}{i}}\}\geq\max\{|J_i^\prime|^{\frac{1}{i}}\}=\h(\p^\prime), \]
which contradicts the assumption. 

Conversely,  let $\p$ be a minimal tuple such that there is no  $\lambda \in \C\setminus\{0\}$ such that $\l^i \in \Z$ and  $\l^i | J_i$ , for all $i =2, \dots , 8$. 
Hence $\p$ cannot be written as   $\left[\l^2J_2^{\prime}, \dots ,\l^8 J_8^{\prime}\right]$, 
where $J_i^\prime \in \Z$ for all $i =2, \dots , 8$.      
Hence,  there is no twist    of $\p$   in $ \in\WP (\Q)$. 
\qed




\def\q{\mathfrak q}

\begin{rem}
Absolute minimal tuples may be not unique over $\overline \Q$. Consider minimal tuples: $[-1 ,-1 ,0 ,-1 ,0 ,0, 0]\in\WP (\Q)$ and $[-1, 1, 0, 1, 0, 0, 0]\in\WP (\Q)$. For both tuples $\h([-1 ,-1 ,0 ,-1 ,0 ,0, 0])=1$ and $\h([-1, 1, 0, 1, 0, 0, 0])=1$. The first tuple is a twist of the second one, with $\lambda=i$: $i\star[-1 ,-1 ,0 ,-1 ,0 ,0, 0]=[-1, 1, 0, 1, 0, 0, 0]$.
\end{rem}

For any absolute minimal tuple $\p$ there may exist other absolute minimal tuples $(-1)\star\p$ or $i\star\p$. 
For any isomorphic class of hyperelliptic curves of genus 3 we want to store only one absolute minimal tuple. For this reason we set a convention that defines which absolute minimal tuples shall be stored.

\begin{description}
	\item[i] For any two corresponding absolute minimal tuples we choose the one with more positive invariants $J_i$.
	\item[ii] For any two corresponding absolute minimal tuples with the same number of positive invariants $J_i$
	we choose the one with greater value $\sum_{i=2}^8i\cdot J_i$.
\end{description}

We will call an absolute minimal tuple that satisfies this convention a \textbf{normalized absolute minimal tuple}. 


\begin{prop}
For every point $\q \in \WP (\Q)$ there is a unique absolute minimal tuple. 
\end{prop}

\proof
Let 
$\p$ and $\p^\prime$ be   absolute  minimal tuples corresponding to $\q \in \WP (\Q)$.  Then there are $\lambda_1, \lambda_2 \in \C$ such that
\[  \q = \lambda_1 \star \p = \lambda_2 \star \p^\prime,\]
and there is no  prime $p$ dividing all coordinates. Hence, 
\[ J_i = \lambda^i J_i^\prime, \; \; i=2, \dots , 8,\]
and $\lambda $ is not divisible by any prime $p$.  Hence, $\lambda $ is a unit in $\Z$ and therefore $\lambda = \pm 1$. If $\lambda =1 $ we have uniqueness, otherwise $\lambda =-1$.  But now the tuple
$[\lambda^2 J_2, \lambda^3 j_3, \dots , \lambda^8 J_8]$ is not an integer tuple any longer and therefore not in $\WP(\Q)$.  

\qed

The above result makes it possible to create a database with all points in the moduli space $\H_3$ with field of moduli $\Q$ and bounded height.  We describe next how this is done, but first let's see a couple of examples.



\begin{exa}
Let $\X$ be the curve with equation 
\[ y^2 = x^7-1.\]
Then the moduli point is 
\[ \p = [0, 0, 0, 0, 0, -395300640, 0] \]
Notice that multiplying this tuple by $\lambda = \frac 1 { (-395300640)^{1/7}}$ gives us the absolute minimal tuple \[ \p = [0, 0, 0, 0, 0, 1, 0].\] 
\end{exa}

Let's see another example.

\begin{exa}
Let $\X$ be the curve with equation 
\[ y^2 = x^8-1.\]
Then the moduli point is 
\[ \p = [-    2^{3} \cdot 5 \cdot 7, 0,   2^{10} \cdot 7^{4}, 0,  2 ^{15} \cdot 7 ^{6}, 0,-  2 ^{19} \cdot 5 \cdot 7^{8}]   \]
Notice that multiplying this tuple by $\lambda = \frac 1  2$ we get a twist  
\[ \p = [-    2  \cdot 5 \cdot 7, 0,   2^6 \cdot 7^{4}, 0,  2^9 \cdot 7^{6}, 0, -  2^{11} \cdot 5 \cdot 7^{8}].  \]
However, we can further reduce its height by multiplying by $\l =  \mathfrak i \frac 1 {\sqrt{2}}$ and we get
\[ \p = [  5 \cdot 7,  0,   2^4 \cdot 7^{4}, 0,  2^6 \cdot 7^{6}, 0, -  2^7 \cdot 5 \cdot 7^{8}].  \]
Furthermore, by multiplying with $\l = \frac 1 {\sqrt{7}}$ we get 
\[ \p = [  5 ,  0,   2^4 \cdot 7^2, 0,  2^6 \cdot 7^3, 0, -  2^7 \cdot 5 \cdot 7^4 ].  \]
which is an absolute minimum tuple. 
\end{exa}

In the next section we show how to create a database of genus 3 hyperelliptic curves which will be ordered by the moduli height. 

\section{A database of genus 3 hyperelliptic curves}

In this part we give a quick computational view on how to create a database of genus $3$ hyperelliptic curves of weighted moduli height $\leq \h$, where $\h$ is a positive integer. We briefly describe each step of our algorithm below. 

\begin{enumerate}
	\item Since we want all the points of height $\leq \h$ that means that we want all ordered tuples $(x_2, \dots , x_8) \in \Z^7$ such that   $|x_i| \leq \h^i$ (tuple with all zeros is not considered). Let the list of all such tuples be denoted by $A$. There are $\prod_{i=2}^8 (2 \h^i +1)-1$ possible tuples in this list. But,  not all   tuples in $A$ are valid tuples for us.  Valid tuples are only the ones that satisfy Eq.~\eqref{shaska} and $J_{14} \neq 0$. 
	
\item Next, for each tuple in $A$ we check if they satisfy the Eq.~\eqref{shaska}:
\begin{equation*}
 J_8^5 + \frac { I_8} {3^4 \cdot 5^3 } J_8^4 + 2 \cdot \frac { I_{16} } {3^8\cdot 5^6}   J_8^3 + \frac {I_{24}} {2 \cdot 3^{12} \cdot 5^6}  J_8^2  + \frac { I_{32}} {3^{16} \cdot 5^{10} }  J_8 + \frac { I_{40}} {2^2 \cdot 3^{20} \cdot 5^{12}}  =0, 
 \end{equation*}
We throw away all tuples which do not satisfy this equation. Let the new set of tuples be denoted by $B$. Elements of set $B$ are points in a weighted projective space $\WP (\Q)$.

\item For each tuple in $B$ we calculate the discriminant $\Delta=J_{14}$. We throw away all tuples with $\Delta=0$ and denote the new set by $C$.

\item Next, for each tuple in set $C$ we compute the corresponding minimal tuple.  
This step is the most expensive one since it requires factorization of each integer coordinate of the tuple into prime factors. Let the new set of tuples be denoted by $D$. 

\item The last step is to compute a normalized absolute minimal tuple for each element in set $D$.
\end{enumerate}

In step (2) we compute invariants $I_8,I_{16},I_{24},I_{32}$ and $I_{40}$ by formulas in \cite{gen-3}.
In step (3) we compute  $\Delta$ as function of $J_2, J_3, \ldots , J_{10}$. 

First, we considered tuples of $\h=1$. 
\begin{enumerate}
	\item A list of all ordered minimal tuples $[x_2, \dots , x_8] \in \Z^7$, such that   $|x_i| \leq 1$ has 2186 elements.
	
\item There are only 34 tuples ($\#B=34$) that satisfy the equation:
\[J_8^5 + a_4 J_8^4 + a_3 J_8^3+ a_2 J_8^2  + a_1 J_8 + a_0 =0.\] 

\item There are 24 tuples with $\Delta\neq0$ ($\#C=24$).

\item Since all elements in $C$ have $\h=1$, $D=C$. 

\item For each minimal tuple in $D$ we compute the corresponding normalized absolute minimal tuples and store them in a database (see Table \ref{nabt1}). 
\end{enumerate}
Notice that  round $ 70\%$ of minimal tuples of height $\h=1$ in a weighted projective space $\WP (\Q)$ correspond to hyperelliptic curves of genus 3. But, only $\approx35\%$ of minimal tuples of hight $\h=1$ in a weighted projective space $\WP (\Q)$ are normalized absolute minimal tuples  that correspond to hyperelliptic curves of genus 3.

\begin{small}
\begin{table}[h]
	\centering
	\caption{Absolute minimal tuples of height 1 (set $C$)}
		\begin{tabular}{c|c|c|c|c|c}
		\#&$\p=[J_2: \ldots : J_8]$ &\#& $\p=[J_2: \ldots : J_8]$ &\#& $\p=[J_2: \ldots : J_8]$\\\hline
   1&-1  -1   0  -1   0   0   0 & 2&-1  -1  0   0   0   0   0 &3 & -1  -1  0  1  0  0  0 \\
	 4&-1  0    0  -1   0   0   0 &5 &-1  0   0   0   0   0   0 &6 & -1  0   0  1  0  0  0 \\
   7&-1  1    0  -1   0   0   0 &8 &-1  1   0   0   0   0   0 &9 & -1  1   0  1  0  0  0 \\
  10& 0 -1    0   -1   0   0   0  &11& 0  -1  0   1   0   0   0 &12&  0  0   0  0  0  -1 0\\
  13& 0  0    0   0   0   1   0 &14& 0  1   0   -1  0   0   0 &15&  0  1   0  1  0  0  0 \\
	16& 1  -1   0  -1   0   0   0 &17& 1  -1  0   0   0   0   0 &18&  1  -1  0  1  0  0  0\\
  19& 1  0    0  -1   0   0   0 &20& 1  0   0   0   0   0   0 &21&  1  0   0  1  0  0  0 \\
	22& 1  1    0  -1   0   0   0 &23& 1  1   0   0   0   0   0 &24 & 1  1   0  1  0  0   0\\
		\end{tabular}
\end{table}
\end{small}

\begin{table}[h]\label{nabt1}
	\centering
	\caption{Normalized absolute minimal tuples of height 1}
		\begin{tabular}{c|c|c|c|c|c}
		\#& $\p=[J_2: \ldots : J_8]$ &\#& $\p=[J_2: \ldots : J_8]$ &\#& $\p=[J_2: \ldots : J_8]$\\\hline
   1& -1  -1   0   1  0   0  0 &2 &-1   0   0    1  0   0   0 &3 &-1  1  0  0  0  0  0 \\
	 4& -1   1   0   1  0   0  0 &5 &0 -1  0  1  0  0  0 & 6& 0   0   0   0  0   1  0\\
   7& 0   1   0    1  0   0   0  & 8&  1 -1  0  1  0  0  0 & 9& 1   0   0   0  0   0  0 \\
	10& 1   0   0    1  0   0   0 &11& 1  1  0  0  0  0  0 &12&1     1     0     1     0     0     0 \\
		\end{tabular}
\end{table}

We considered tuples of $\h\leq 1.5$. 
\begin{enumerate}
	\item A list of all ordered minimal tuples $[x_2, \dots , x_8] \in \Z^7$, such that   $|x_i| \leq (1.5)^i$ has 237092624 elements.
	
\item There are only 748 tuples ($\#B=748$) that satisfy the equation:
\[J_8^5 + a_4 J_8^4 + a_3 J_8^3+ a_2 J_8^2  + a_1 J_8 + a_0 =0.\] 

\item There are 544 tuples with $\Delta\neq0$ ($\#C=544$).

\item All elements in $C$ are minimal tuples so $D=C$. 

\item For each minimal tuple in $D$ we compute the corresponding normalized absolute minimal tuple and store them in a database (see Table \ref{nabt2}). There is $246+12=258$ normalized absolute minimal tuples of hight $0<\h\leq1.5$.
\end{enumerate}
Notice that  round $ 73.(73)\%$ of minimal tuples of height $1\leq\h\leq1.5$ in a weighted projective space $\WP (\Q)$ correspond to hyperelliptic curves of genus 3. But, only $\approx34\%$ of minimal tuples of hight $1\leq\h\leq1.5$ in a weighted projective space $\WP (\Q)$ are normalized absolute minimal tuples  that correspond to hyperelliptic curves of genus 3.

\begin{table}[h]\label{nabt2}
	\centering
\caption{Normalized absolute minimal tuples of height $1<\h\leq1.5$}
		\begin{tabular}{c|c|c|c|c|c}
		\#& $\p=[J_2: \ldots : J_8]$ &\#& $\p=[J_2: \ldots : J_8]$ &\#& $\p=[J_2: \ldots : J_8]$\\\hline
   1& -2	-3	0	1	0	0	0 & 2 &-2	-3	0	2	0	0	0& 3& -2	-3	0	3	0	0	0\\
4&-2	-3	0	4	0	0	0& 5& -2	-3	0	5	0	0	0& 6& -2	-3	0	6	0	0	0\\
7&-2	-3	0	7	0	0	0& 8& -2	-2	0	1	0	0	0& 9& -2	-2	0	2	0	0	0\\
10&-2	-2	0	3	0	0	0& 11& -2	-2	0	4	0	0	0& 12&-2	-2	0	5	0	0	0\\
13&-2	-2	0	6	0	0	0& 14& -2	-2	0	7	0	0	0 & 15&-2	-1	0	1	0	0	0\\
16&-2	-1	0	2	0	0	0& 17& -2	-1	0	3	0	0	0& 18& -2	-1	0	4	0	0	0\\
19&-2	-1	0	5	0	0	0& 20& -2	-1	0	6	0	0	0& 21& -2	-1	0	7	0	0	0\\
22&-2	0	0	1	0	0	0&23& -2	0	0	2	0	0	0& 24&-2	0	0	3	0	0	0\\
25&-2	0	0	4	0	0	0& 26& -2	0	0	5	0	0	0& 27& -2	0	0	6	0	0	0\\
28&-2	0	0	7	0	0	0& 29& -2	1	0	0	0	0	0& 30& -2	1	0	1	0	0	0\\
31&-2	1	0	2	0	0	0& 32&-2	1	0	3	0	0	0 &33&-2	1	0	4	0	0	0\\
34&-2	1	0	5	0	0	0& 35& -2	1	0	6	0	0	0& 36&-2	1	0	7	0	0	0\\
37&-2	2	0	0	0	0	0& 38& -2	2	0	1	0	0	0 &39&-2	2	0	2	0	0	0\\
40&-2	2	0	3	0	0	0& 41& -2	2	0	4	0	0	0& 42& -2	2	0	5	0	0	0\\
43&-2	2	0	6	0	0	0&44&-2	2	0	7	0	0	0 &45& -2	3	0	0	0	0	0\\
46&-2	3	0	1	0	0	0 &47& -2	3	0	2	0	0	0 &48&-2	3	0	3	0	0	0\\
49&-2	3	0	4	0	0	0& 50& -2	3	0	5	0	0	0& 51&-2	3	0	6	0	0	0\\
52&-2	3	0	7	0	0	0 &53&-1	-3	0	1	0	0	0 &54&-1	-3	0	2	0	0	0\\
55&-1	-3	0	3	0	0	0& 56&-1	-3	0	4	0	0	0& 57&-1	-3	0	5	0	0	0\\
58&-1	-3	0	6	0	0	0& 59&-1	-3	0	7	0	0	0 &60&-1	-2	0	1	0	0	0\\
61& -1	-2	0	2	0	0	0&62&-1	-2	0	3	0	0	0& 63&-1	-2	0	4	0	0	0 \\
64&-1	-2	0	5	0	0	0&65& -1	-2	0	6	0	0	0&66&-1	-2	0	7	0	0	0\\
67&-1	-1	0	2	0	0	0 &68&-1	-1	0	3	0	0	0&69&-1	-1	0	4	0	0	0\\
70& -1	-1	0	5	0	0	0 &71&-1	-1	0	6	0	0	0& 72&-1	-1	0	7	0	0	0\\
73&  -1	0	0	2	0	0	0&74& -1	0	0	3	0	0	0&75&-1	0	0	4	0	0	0\\
76& -1	0	0	5	0	0	0&77& -1	0	0	6	0	0	0&78&-1	0	0	7	0	0	0\\
79& -1	1	0	2	0	0	0&80& -1	1	0	3	0	0	0&81&-1	1	0	4	0	0	0\\
82& -1	1	0	5	0	0	0&83&-1	1	0	6	0	0	0& 84&-1	1	0	7	0	0	0 \\
85&-1	2	0	0	0	0	0 &86&-1	2	0	1	0	0	0&87&-1	2	0	2	0	0	0 \\
88&-1	2	0	3	0	0	0&89& -1	2	0	4	0	0	0&90&-1	2	0	5	0	0	0\\
91& -1	2	0	6	0	0	0&92& -1	2	0	7	0	0	0 &93&-1	3	0	0	0	0	0\\
94& -1	3	0	1	0	0	0 &95&-1	3	0	2	0	0	0 &96&-1	3	0	3	0	0	0\\
97&  -1	3	0	4	0	0	0 &98&-1	3	0	5	0	0	0&99&-1	3	0	6	0	0	0  \\
100& -1	3	0	7	0	0	0&101& 0	-3	0	1	0	0	0&102&0	-3	0	2	0	0	0\\
103& 0	-3	0	3	0	0	0 &104&0	-3	0	4	0	0	0&105&0	-3	0	5	0	0	0\\
106&  0	-3	0	6	0	0	0&107& 0	-3	0	7	0	0	0&108&0	-2	0	1	0	0	0 \\
109&0	-2	0	2	0	0	0&110& 0	-2	0	3	0	0	0&111& 0	-2	0	4	0	0	0 \\
112&0	-2	0	5	0	0	0 &113& 0	-2	0	6	0	0	0&114& 0	-2	0	7	0	0	0\\
115& 0	-1	0	2	0	0	0&116& 0	-1	0	3	0	0	0 &117&0	-1	0	4	0	0	0 \\
118&0	-1	0	5	0	0	0 &119& 0	-1	0	6	0	0	0&120&0	-1	0	7	0	0	0\\
121& 0	0	-3	0	-10	0	-5&122&  0	0	-3	0	10	0	-5 &123&0	0	3	0	-8	0	-5 \\
124& 0	0	3	0	-2	0	-5&125&  0	0	3	0	2	0	-5 &126& 0	0	3	0	8	0	-5\\
127&  0	1	0	2	0	0	0 &128&0	1	0	3	0	0	0&129& 0	1	0	4	0	0	0 \\
130& 0	1	0	5	0	0	0&131& 0	1	0	6	0	0	0&132& 0	1	0	7	0	0	0\\
133&0	2	0	1	0	0	0 &134& 0	2	0	2	0	0	0 &135&0	2	0	3	0	0	0 \\
136& 0	2	0	4	0	0	0 &137& 0	2	0	5	0	0	0 &138&0	2	0	6	0	0	0  \\
139& 0	2	0	7	0	0	0&140& 0	3	0	1	0	0	0 &141& 0	3	0	2	0	0	0\\
		\end{tabular}
\end{table}
\begin{table}[h]
	\centering
	\caption{Normalized absolute minimal tuples of height $1<\h\leq1.5$}
		\begin{tabular}{c|c|c|c|c|c}
		\#& $\p=[J_2: \ldots : J_8]$ &\#& $\p=[J_2: \ldots : J_8]$ &\#& $\p=[J_2: \ldots : J_8]$\\\hline
142& 0	3	0	3	0	0	0&143&  0	3	0	4	0	0	0 &144&0	3	0	5	0	0	0 \\
145&0	3	0	6	0	0	0&146&  0	3	0	7	0	0	0 &147&1	-3	0	1	0	0	0\\
148& 1	-3	0	2	0	0	0&149& 1	-3	0	3	0	0	0 &150&1	-3	0	4	0	0	0\\
151&1	-3	0	5	0	0	0&152 &  1	-3	0	6	0	0	0 &153&1	-3	0	7	0	0	0\\
154& 1	-2	0	1	0	0	0 &155& 1	-2	0	2	0	0	0 &156&1	-2	0	3	0	0	0\\
157&  1	-2	0	4	0	0	0&158& 1	-2	0	5	0	0	0 &159&1	-2	0	6	0	0	0 \\
160&1	-2	0	7	0	0	0&161& 1	-1	0	2	0	0	0 &162&1	-1	0	3	0	0	0 \\
163& 1	-1	0	4	0	0	0&164& 1	-1	0	5	0	0	0&165&1	-1	0	6	0	0	0\\
166& 1	-1	0	7	0	0	0&167& 1	0	0	2	0	0	0 &168&1	0	0	3	0	0	0 \\
169&   1	0	0	4	0	0	0 &170&1	0	0	5	0	0	0&171&1	0	0	6	0	0	0\\
172& 1	0	0	7	0	0	0&173&  1	1	0	2	0	0	0 &174&1	1	0	3	0	0	0\\
175&  1	1	0	4	0	0	0&176 &  1	1	0	5	0	0	0 &177&1	1	0	6	0	0	0\\
178& 1	1	0	7	0	0	0 &179&1	2	0	0	0	0	0 &180& 1	2	0	1	0	0	0\\
181& 1	2	0	2	0	0	0  &182&1	2	0	3	0	0	0&183&1	2	0	4	0	0	0\\
184 &1	2	0	5	0	0	0 &185&  1	2	0	6	0	0	0&186& 1	2	0	7	0	0	0\\
187&1	3	0	0	0	0	0  &188& 1	3	0	1	0	0	0&189&1	3	0	2	0	0	0  \\
190&1	3	0	3	0	0	0 &191&  1	3	0	4	0	0	0&192&1	3	0	5	0	0	0\\
193&  1	3	0	6	0	0	0 &194& 1	3	0	7	0	0	0&195& 2	-3	0	1	0	0	0  \\
196& 2	-3	0	2	0	0	0&197&2	-3	0	3	0	0	0&198&2	-3	0	4	0	0	0 \\
199& 2	-3	0	5	0	0	0&200&  2	-3	0	6	0	0	0 &201&2	-3	0	7	0	0	0\\
202& 2	-2	0	1	0	0	0&203& 2	-2	0	2	0	0	0 &204&2	-2	0	3	0	0	0 \\
205& 2	-2	0	4	0	0	0 &206& 2	-2	0	5	0	0	0 &207& 2	-2	0	6	0	0	0\\
208& 2	-2	0	7	0	0	0 &209& 2	-1	0	1	0	0	0 &210&2	-1	0	2	0	0	0\\
211& 2	-1	0	3	0	0	0&212& 2	-1	0	4	0	0	0 &213& 2	-1	0	5	0	0	0 \\
214&  2	-1	0	6	0	0	0&215& 2	-1	0	7	0	0	0&216& 2	0	0	1	0	0	0\\
217& 2	0	0	2	0	0	0 &218& 2	0	0	3	0	0	0 &219&2	0	0	4	0	0	0\\
220& 2	0	0	5	0	0	0 &221& 2	0	0	6	0	0	0 &222&2	0	0	7	0	0	0  \\
223&2	1	0	0	0	0	0&224& 2	1	0	1	0	0	0 &225&2	1	0	2	0	0	0  \\
226& 2	1	0	3	0	0	0 &227& 2	1	0	4	0	0	0 &228&2	1	0	5	0	0	0 \\
229& 2	1	0	6	0	0	0&230 &  2	1	0	7	0	0	0&231& 2	2	0	0	0	0	0\\
232& 2	2	0	1	0	0	0 &233&  2	2	0	2	0	0	0&234&2	2	0	3	0	0	0\\
235& 2	2	0	4	0	0	0 &236& 2	2	0	5	0	0	0 &237& 2	2	0	6	0	0	0 \\
238& 2	2	0	7	0	0	0 &239&2	3	0	0	0	0	0 &240& 2	3	0	1	0	0	0 \\
241& 2	3	0	2	0	0	0&242&2	3	0	3	0	0	0 &243& 2	3	0	4	0	0	0 \\
244& 2	3	0	5	0	0	0&245& 2	3	0	6	0	0	0&246&2	3	0	7	0	0	0\\
		\end{tabular}
\end{table}


\clearpage

\nocite{*}

\bibliographystyle{amsacm}

\bibliography{ref}{}
\end{document}